\documentclass[12pt,a4paper]{amsart}
\usepackage{amssymb, amsmath, amsthm}
\usepackage[all]{xy}

\theoremstyle{plain}
\newtheorem{thm}{Theorem}
\newtheorem{cor}{Corollary}
\newtheorem{prop}{Proposition}
\newtheorem{lem}{Lemma}

\theoremstyle{definition}

\theoremstyle{remark}
\newtheorem{rem}[subsubsection]{Remark}
\newtheorem{exam}[subsubsection]{Example}
\newtheorem{prob}{Problem}

\DeclareMathOperator{\dist}{dist}
\DeclareMathOperator{\Log}{Log}

\newcommand{\imax}{\underline{m}}

\title[Quasianalytic local rings]{On quasianalytic local rings}

\author{Vincent Thilliez}

\address{Math\'ematiques - B\^atiment M2\\
Universit\'e des Sciences et Technologies de Lille\\
F-59655 Villeneuve d'Ascq Cedex, France}

\email{thilliez@math.univ-lille1.fr}


\begin{document}

\begin{abstract} This expository article is devoted to the local theory of ultradifferentiable classes of functions, with a special emphasis on the quasianalytic case. Although quasianalytic classes are well-known in harmonic analysis since several decades, their study from the viewpoint of differential analysis and analytic geometry has begun much more recently and, to some extent, has earned them a new interest. Therefore, we focus on contemporary questions closely related to topics in local algebra. We study, in particular, Weierstrass division problems and the role of hyperbolicity, together with properties of ideals of quasianalytic germs. Incidentally, we also present a simplified proof of Carleman's theorem on the non-surjectivity of the Borel map in the quasianalytic case.
\end{abstract}

\maketitle

\section{Ultradifferentiable function germs}

\subsection{Historical background}\label{backgnd}
At the end of the nineteenth century \cite{Borel1, Borel2}, Borel produced the first non-trivial examples of sets $ E $ of infinitely differentiable functions on the real line, containing nowhere analytic functions, and such that any element $ f $ in $ E $ satisfies the implication 
\begin{equation}\label{QAhisto}
(f^{(j)}(0)=0,\ j=0,1,\ldots )\Longrightarrow (f=0). 
\end{equation} 
Borel's examples were typically given by restrictions to the real line of series of rational functions 
\begin{equation}\label{Borelex}
f(z)=\sum_{\nu=1}^{+\infty}\frac{A_\nu}{z-z_\nu},
\end{equation}
where the poles $ z_\nu $ belong to $ \mathbb{C}\setminus\mathbb{R} $ and accumulate near the real axis, whereas the coefficients $ A_\nu $ tend sufficiently fast to $ 0 $.  

\begin{exam} Here is a simplified realization of such examples (the result is weaker since non-analyticity is merely achieved at $ 0 $). Put $ z_\nu=i/\nu $ and consider the set $ E $ of functions $ f $ defined by \eqref{Borelex}, with the requirement $\varlimsup \vert A_\nu\vert^{1/\nu}<1 $. Being given such a sequence $ (A_\nu)_{\nu\geq 1} $, put $ \Phi(\zeta)=\sum_{\nu=1}^{+\infty}\nu A_\nu \zeta^\nu $. Observe that $ \Phi $ is holomorphic in a disc of radius strictly greater than $ 1 $, hence
$ \sum_{\nu= 1}^{+\infty}\nu^{j+1}\vert A_\nu\vert<\infty $ for any integer $ j $. It is then easy to see that $ f $ is a $ C^\infty $ function on the real line and that its derivatives at $ 0 $ are linear combinations of the derivatives of $ \Phi $ at $ 1 $. These linear combinations are given by an invertible triangular matrix, hence property \eqref{QAhisto} follows. One can also check that the series defining $ f $ converges to a holomorphic function on $ \mathbb{C}\setminus \Delta $, where $ \Delta $ is the half-line $ \{it : t\geq 0\} $. If we assume $ A_\nu\neq 0 $ for every $ \nu $, this function cannot be extended to a neighborhood of $ 0 $ in $ \mathbb{C} $, hence $ f $ is not analytic in a neighborhood of $ 0 $ in $ \mathbb{R} $.
\end{exam}

About a decade after Borel's discovery, Ha\-da\-mard was led to a decisive idea by considerations of PDE theory. Indeed, the work of Holmgren on the heat equation had already revealed that the solutions of certain partial differential equations are natural elements of classes of functions between analytic and $ C^\infty $ ones, defined by bounds on their successive derivatives (they are now well-known as \textit{Gevrey classes}, by reference to \cite{Gevrey}). Hadamard, in a communication to the \textit{Soci\'et\'e Math\'ematique de France} \cite{Hadamard}, asked whether the implication \eqref{QAhisto} could also be characterized in terms of a growth condition on derivatives. As we shall see, the answer is essentially affirmative and serves as a starting point to a number of results.

\subsection{General framework}
We use the following notation: for any multi-index $ J=(j_1,\ldots,j_n) $ of $ \mathbb{N}^n $, we denote the length $ j_1+\cdots+j_n $ of $ J $ by the corresponding lower case letter $ j $. We put 
$ D^J=\partial^j/\partial x_1^{j_1}\cdots\partial x_n^{j_n} $, $ J!=j_1!\cdots j_n! $ and $ x^J=x_1^{j_1}\cdots x_n^{j_n} $.

Frow now on, unless explicitly mentioned, we shall work from a local viewpoint. The customary identification between germs of sets or mappings and their representatives will always be used, unless there is a danger of confusion.

Denote by $ \mathcal{O}_n $ (resp. $ \mathcal{E}_n $) the ring of real-analytic (resp. infinitely differentiable), complex-valued, function germs at the origin of $ \mathbb{R}^n $, and by $ \mathcal{F}_n $ the ring of formal power series in $ n $ variables, with complex coefficients. The map $ T_0 : \mathcal{E}_n\longrightarrow \mathcal{F}_n $ defined by
\begin{equation*}
T_0f=\sum_{J\in \mathbb{N}^n}\frac{D^Jf(0)}{J!}x^J
\end{equation*}
will be called the \emph{Borel map}.

Now, let $ M=(M_j)_{j\geq 0} $ be an increasing sequence of real numbers, with $ M_0=1 $. Denote by $ \mathcal{E}_n(M) $ the set of elements $ f $ of $ \mathcal{E}_n $ for which there exist a neighborhood $ U $ of $ 0 $ and positive constants $ C $ and $ \sigma $ such that
\begin{equation}\label{carleman}
\big\vert D^Jf(x)\big\vert \leq C\sigma^jj!M_j\ \textrm{ for any } J\in \mathbb{N}^n\textrm{ and } x\in U.
\end{equation}
Here, $ C\sigma^j j! $ appears as ``the analytic part'' of the estimate, whereas $ M_j $ can be considered as a way to allow a defect of analyticity. We clearly have 
\begin{equation*}
\mathcal{O}_n\subseteq\mathcal{E}_n(M)\subseteq\mathcal{E}_n.
\end{equation*}
In Corollary \ref{nonana} hereafter, we shall characterize the case of equality in the first inclusion. The second inclusion is, obviously, always strict. 

In the same spirit, denote now by $ \mathcal{F}_n(M) $ the set of elements $ F=\sum_{J\in\mathbb{N}^n}F_Jx^J $ of $\mathcal{F}_n $ for which there exist positive constants $ C $ and $ \sigma $ such that
\begin{equation*}
\vert F_J\vert\leq C\sigma^j M_j\ \textrm{ for any } J\in \mathbb{N}^n.
\end{equation*} 
The Borel map then obviously satisfies 
\begin{equation*}
T_0\mathcal{E}_n(M)\subseteq \mathcal{F}_n(M).
\end{equation*}

\subsection{Basic structure}\label{basicstru}
One cannot hope to get much more information on the sets $ \mathcal{E}_n(M) $ and $ \mathcal{F}_n(M) $ without an assumption on the sequence $ M $. Frow now on, we shall \emph{always} make the following assumption:
\begin{equation}\label{logconv}
\textrm{ the sequence } M \textrm{ is logarithmically convex}.
\end{equation}
This amounts to saying that $ M_{j+1}/M_j $ increases. Taking into account the value $ M_0=1 $, it is easy to derive that $ (M_j)^{1/j} $ also increases, and 
\begin{equation}\label{prod}
M_jM_k\leq M_{j+k}\ \textrm{ for any } (j,k)\in\mathbb{N}^2.
\end{equation}
Assumption \eqref{logconv} implies several important properties. 
First, the set $ \mathcal{E}_n(M) $ is a ring, as can be derived from \eqref{prod} by means of the Leibniz formula\footnote{In fact, it suffices to assume the logarithmic convexity of the sequence $ \overline{M} $ defined by $ \overline{M}_j= j!M_j $, but we shall not use this refinement here.}. Second, stability under composition holds in the following sense: for any $ f=(f_1,\ldots,f_p) $ in $ (\mathcal{E}_n(M))^p $, with $ f(0)=0 $, and any $ g $ in $ \mathcal{E}_p(M) $, the composite $ g\circ f $ belongs to $ \mathcal{E}_n(M) $. In this general form, the result is due to Roumieu \cite{Roumieu}, although particular cases are much older \cite{Cartan1, Gevrey}. As a consequence, we obtain the following fact.

\begin{prop}
The set $ \mathcal{E}_n(M) $ is a local ring with maximal ideal $ \imax_M=\{h\in \mathcal{E}_n(M) : h(0)=0\} $.
\end{prop}

\begin{proof}
If an element $ h $ of $ \mathcal{E}_n(M) $ is such that $ h(0)\neq 0 $, then it is invertible in $ \mathcal{E}_n(M) $, as can be seen by composition of $ f(x)=h(x)-h(0) $ and $ g(t)=(h(0)+t)^{-1} $. Indeed, $ g$ belongs to $ \mathcal{O}_1 $, hence to $ \mathcal{E}_1(M) $.
\end{proof}

\begin{rem}
Similar techniques can be used in the formal case, hence $ \mathcal{F}_n(M) $ is also a local ring. 
\end{rem}

Among other consequences of \eqref{logconv}, we mention the implicit function theorem for $ \mathcal{E}_n(M) $; see \cite{Komatsu2}. It should be remarked that stability under composition and the implicit function theorem hold, in fact, under slightly weaker assumptions, which are not the same for both results. In the same way, \eqref{logconv} can also be weakened in some of the statements hereafter, whereas other ones actually require the assumption. We do not consider these issues, since we favor a unified treatment.

\subsection{Equivalence and stability properties}
Here we discuss a characterization of the inclusion $ \mathcal{E}_n(M)\subseteq\mathcal{E}_n(N) $ by a condition on the sequences $ M $ and $ N $. We write $ M\prec N $ if there exists a constant $ C $ such that $ M_j\leq C^j N_j $ for any $ j\geq 0 $, or, equivalently, if $ \sup_{j\geq 1}(M_j/N_j)^{1/j}<\infty $. Clearly, the condition $ M\prec N $ implies $\mathcal{E}_n(M)\subseteq \mathcal{E}_n(N)$ and $\mathcal{F}_n(M)\subseteq \mathcal{F}_n(N)$. The converse implication is clear in the formal case: if the element $ F=\sum_{j=0}^{+\infty}M_jx_n^j $ of $ \mathcal{F}_n(M) $ belongs to $ \mathcal{F}_n(N) $, we immediately derive $ M\prec N $. This argument does not work in the case of germs, since there is no reason why $ F $ should be the Taylor series of some element $ f$ of $ \mathcal{E}_n(M) $ (as a matter of fact, we shall see in Section \ref{propBor} that $ T_0 $ is generally not surjective). Fortunately, $ \mathcal{E}_n(M) $ contains functions with sufficiently large derivatives\footnote{In particular, such functions do not belong to any smaller local ring of the same type.}, as shown by the following variant of classical results of Cartan \cite{Cartan1} and Mandelbrojt \cite{Cartan2}.

\begin{thm}\label{cartan} 
There exists an element $ \theta $ in $ \mathcal{E}_1(M) $ such that 
$ \big\vert\theta^{(j)}(0)\big\vert\geq j!M_j $ for any integer $ j\geq 0 $. 
\end{thm}

\begin{proof} 
We use a construction originating in the classical work of Bang \cite{Bang}. Put $ \overline{M}_j=j!M_j $ and $ m_j=\overline{M}_{j+1}/\overline{M}_j $. The logarithmic convexity of $ M $ implies that $ (\overline{M}_j)_{j\geq 1} $ is also logarithmically convex, hence $ (m_j)_{j\geq 1} $ increases. Discussing separately the cases $ j\leq k $ and $ j>k $, it is easy to obtain the estimate
\begin{equation*}
\left(\frac{1}{m_k}\right)^{k-j}\leq\frac{\overline{M}_j}{\overline{M}_k}\ \textrm{ for any } (j,k)\in \mathbb{N}^2.
\end{equation*}
It is then not difficult to check that the function $ \theta $ defined by 
\begin{equation*}
\theta(x)=\sum_{k=0}^{+\infty}\frac{\overline{M}_k}{(2m_k)^k}\exp\left(2im_k x\right)
\end{equation*}
has all the required properties.
\end{proof}

Theorem \ref{cartan} implies that we have $ \mathcal{E}_n(M)\subseteq \mathcal{E}_n(N) $ if and only if $ M\prec N $. Replacing $ N_j $ by $ 1 $, we get the following result.

\begin{cor}\label{nonana}
We have $ \mathcal{O}_n =\mathcal{E}_n(M) $ if and only if $ \sup_{j\geq 1}(M_j)^{1/j}<\infty $.
\end{cor}

\begin{rem}\label{nonan}
As pointed out in \S \ref{basicstru}, the logarithmic convexity of $ M $ implies that  $ (M_j)^{1/j} $ increases. The strict inclusion $  \mathcal{O}_n \subsetneq\mathcal{E}_n(M) $ is therefore equivalent to the condition $ \lim_{j\to \infty}(M_j)^{1/j}=\infty $. 
\end{rem}

Remark now that the first order derivatives of the elements of $ \mathcal{E}_n(M) $ belong to $ \mathcal{E}_n(M') $ with $ M'_j=M_{j+1} $. Thus, replacing $ M_j $ by $ M'_j $ and $ N_j $ by $ M_j $, we get another important corollary. 

\begin{cor}\label{stabder}
The ring $ \mathcal{E}_n(M) $ is stable under derivation if and only if the sequence $ M $ satisfies
$ \sup_{j\geq 1}(M_{j+1}/M_j)^{1/j}<\infty $. 
\end{cor}

\begin{rem} In this case, it is easy to check that the maximal ideal $ \imax_M $ is generated by the coordinate functions $ x_1,\ldots,x_n $.
\end{rem}

\section{On the Borel map}\label{propBor}

\subsection{Injectivity vs. Surjectivity}
The local ring $ \mathcal{E}_n(M) $ is said to be \emph{quasianalytic} if all its elements satisfy \eqref{QAhisto}. This amounts to saying that the Borel map $ T_0 $ is injective on  $\mathcal{E}_n(M)$. Quasianalyticity is characterized by the famous Denjoy-Carleman theorem, which therefore answers to the question of Hadamard mentioned in \S \ref{backgnd}. The theorem was proved in several steps between 1921 and 1923 \cite{Carleman1, Carleman12, Denjoy}. In our language, it can be stated as follows.

\begin{thm}\label{DCthm}
The local ring $ \mathcal{E}_n(M) $ is quasianalytic if and only if 
\begin{equation}\label{QAcond}
\sum_{j=0}^{+\infty}\frac{M_j}{(j+1)M_{j+1}}=\infty.
\end{equation}
\end{thm}

Contemporary proofs of the theorem can be found in many places; see for instance \cite{Hormander, Katznelson, Koosis, Rudin} (remark that it is easy to reduce the problem to the case $ n=1 $ treated in all these references). The following result, which is essentially a restatement of another theorem of Carleman \cite{Carleman11, Carleman2}, provides an important additional information.

\begin{thm}\label{nonsurj} 
Assume that $ \mathcal{E}_n(M) $ is quasianalytic and that $ \mathcal{O}_n\subsetneq\mathcal{E}_n(M) $. Then the map $ T_0 : \mathcal{E}_n(M)\longrightarrow \mathcal{F}_n(M) $ is not surjective.
\end{thm}

Carleman's original proof of Theorem \ref{nonsurj} relies on a somewhat delicate variational argument involving integral equations. In \S \ref{revisit} below, we shall present a more elementary proof, using Hilbert space techniques to bypass this argument. 

In view of Theorem \ref{nonsurj}, it is natural to ask for a characterization of surjectivity. In other words, when is there a Borel-type extension theorem from $ \mathcal{F}_n(M) $ to $ \mathcal{E}_n(M) $? As we shall see, non-quasianalyticity is not enough. From the early Sixties, the problem was studied a lot, until Petzsche eventually provided a complete solution\footnote{In fact, Petzsche's paper does more, since it also deals with continuous extension operators, but this is another story. The article is written for functions of one variable exclusively, but similar constructions can be carried out in several variables: see, for instance, \cite{Beaugendre}, Section 5.2.} in 1988; see \cite{Petzsche} and the references therein. Petzsche's result can be stated as follows.

\begin{thm}\label{T0surj}
Assume that $ \mathcal{O}_n\subsetneq \mathcal{E}_n(M) $. Then the map $ T_0 : \mathcal{E}_n(M)\longrightarrow \mathcal{F}_n(M) $ is surjective if and only if there exists a constant $ C $ such that
\begin{equation}\label{SNQA}
\sum_{j= k}^{+\infty}\frac{M_j}{(j+1)M_{j+1}}\leq C\frac{M_k}{M_{k+1}}\ \textrm{ for any integer } k\geq 0.
\end{equation}
\end{thm}

Condition \eqref{SNQA} is known as the \emph{strong 
non-quasianalyticity} condition. Under the (obviously weaker) non-quasianalyticity assumption, it is only possible to recover extensions with a loss of smoothness, that is, from $ \mathcal{F}_n(M) $ to $ \mathcal{E}_n(N) $ for some other (larger) sequence $ N $ suitably related to $ M $. See, for instance, \cite{Schmets} and the references therein.
 
\begin{exam}  
Let $\alpha $ be a real number, with $ \alpha>0 $. Put $ M_j=(\Log (j+e))^{\alpha j} $. Then $ \mathcal{E}_n(M) $ is quasianalytic for $ 0<\alpha\leq 1 $ and non-quasianalytic (but not strongly) for $ \alpha>1 $.  
\end{exam}

\begin{exam}  
Let $\alpha $ be a real number, with $ \alpha>0 $. Put $ M_j=(j!)^\alpha $. Then $ \mathcal{E}_n(M) $ is strongly non-quasianalytic. This is the Gevrey $ G^{1+\alpha} $ regularity well-known in PDE theory. 
\end{exam}

\begin{exam}  
Let $ q $ be a real number, with $ q>1 $. Put $ M_j=q^{j^2} $. Then $ \mathcal{E}_n(M) $ is strongly non-quasianalytic. This is the $ q$-Gevrey regularity arising in the study of difference equations. 
\end{exam}

\subsection{Intricacies of the quasianalytic setting}
Theorem \ref{nonsurj} reveals an inherent difficulty of the quasianalytic case. In fact, one currently knows no concrete way to recognize the elements of $ T_0\mathcal{E}_n(M) $ among those of $ \mathcal{F}_n(M) $. Thus, computations on the Taylor series of germs do not say much on the germs themselves.  Another striking phenomenon appears in the following theorem of Mandelbrojt \cite{Mandelbrojt}.

\begin{thm}\label{mandelsum} 
Let $ f $ be an element of $ \mathcal{E}_n $. Then there exist quasianalytic local rings $ \mathcal{E}_n(M) $ and $ \mathcal{E}_n(N) $ (both depending on $ f $) such that 
$ f=g+h $ with $ g\in \mathcal{E}_n(M) $ and $ h\in \mathcal{E}_n(N) $.
\end{thm}

Remark that Mandelbrojt's original proof does not provide the logarithmic convexity property \eqref{logconv} that we require. A proof taking this additional feature into account can be found in \cite{RSW}, where Theorem \ref{mandelsum} finds unexpected application in model theory\footnote{Precisely, it is used to show that there is no largest o-minimal expansion of the real field.}. 

\subsection{A proof of Theorem \ref{nonsurj}}\label{revisit}
As announced, we shall now present an elementary proof of Carleman's theorem on the failure of surjectivity for the Borel map in the quasianalytic setting. We need first to define some function spaces that will also be used in the next sections. 

\subsubsection*{Function spaces}
For any integer $ \nu\geq 1 $ and any real $ \sigma>0 $, we put $ I_\nu=]-1/\nu,1/\nu[ $ and 
$ \mathcal{E}_{n,\sigma}(M,\nu)=\{f\in C^\infty(I_\nu^n) : \Vert f\Vert_{\sigma,\nu}<\infty \} $ with
\begin{equation*}
\Vert f\Vert_{\sigma,\nu}=\sup_{J\in\mathbb{N}^n,\ x\in I_\nu^n}\frac{\vert D^Jf(x)\vert}{\sigma^j j!M_j}.
\end{equation*}
It is easy to check that $ \mathcal{E}_{n,\sigma}(M,\nu) $ is a Banach space. 

\begin{rem}\label{inclucomp}
For $ \sigma<\sigma' $, the canonical injection $ \mathcal{E}_{n,\sigma}(M,\nu)\hookrightarrow \mathcal{E}_{n,\sigma'}(M,\nu) $ is known to be a compact map: this is a consequence of Ascoli's theorem (see \cite{Komatsu}, or \cite{Bourbaki} in the Gevrey case).
\end{rem}

\subsubsection*{Scheme of the proof} 
Clearly, it is enough to deal with the case $n=1 $. The proof then consists of two parts. In the first part, a representation formula for quasianalytic germs is established; this is where we replace the methods of \cite{Carleman2}, Chapter VII, by more elementary ones. The second part follows more closely Carleman's original argument.

\subsubsection*{First part}
Put $ I=I_1=]-1,1[ $ and denote by $ C^\infty(\overline{I}) $ the space of functions $ u $ which are smooth in $ I $ and whose derivatives $ u^{(j)} $, at any order $j \geq 0 $, extend continuously to $ \overline{I} $. For any function $ v $ continuous in $ I $, put $ \Vert v\Vert_{L^2(I)}=\left(\int_I\vert v\vert^2\right)^{1/2} $ and $ \Vert v\Vert_{L^\infty(I)}=\sup_{x\in I}\vert v(x)\vert $. Then, for any element $ u $ of $ C^\infty(\overline{I}) $ and any integer $ j\geq 0 $, we have the elementary inequalities\footnote{Indeed, given $ v $ in $ C^\infty(\overline{I}) $, there exists $ c $ in $ I $ such that $ \vert v(c)\vert=\frac{1}{2}\int_I\vert v\vert  $. For any $ x $ in $ I $, we have $ \vert v(x)\vert\leq \vert v(c)\vert+\left\vert \int_c^x\vert v'\vert \right\vert \leq \frac{1}{2}\int_I\vert v\vert +\int_I\vert v'\vert $, hence $ \vert v(x)\vert\leq \frac{\sqrt{2}}{2}\left(\int_I\vert v\vert^2\right)^{1/2}+\sqrt{2}
\left(\int_I\vert v'\vert^2\right)^{1/2} $ by the Cauchy-Schwarz inequality. Taking $ v=u^{(j)} $, we get the right-hand side of \eqref{sobo}. The left-hand side is trivial.} 
\begin{equation}\label{sobo}
\frac{1}{\sqrt{2}}\big\Vert u^{(j)}\big\Vert_{L^2(I)}\leq \big\Vert u^{(j)}\big\Vert_{L^\infty(I)}\leq \sqrt{2}\left( \big\Vert u^{(j)}\big\Vert_{L^2(I)}+\big\Vert u^{(j+1)}\big\Vert_{L^2(I)}\right).
\end{equation}
Now, put    
\begin{equation*}
\Vert u\Vert^2=\sum_{j=0}^{+\infty}\left(j!M_j\right)^{-2}\big\Vert u^{(j)}\big\Vert_{L^2(I)}^2
\end{equation*}
and denote by $ \mathcal{H}_M $ the space of those $ u $ such that $ \Vert u\Vert^2<\infty $.
It is easy to check that $ \mathcal{H}_M $ is a Hilbert space for the norm $ \Vert\cdot\Vert $; the associated scalar product will be denoted by $ \langle\cdot\vert\cdot\rangle $. Put $ M'_j=M_{j+1} $. Using \eqref{sobo}, we obtain, for any real number $ \eta\in ]0,1[ $, the topological inclusion
\begin{equation}\label{inclusob}
\mathcal{E}_{1,1-\eta}(M,1)\subseteq \mathcal{H}_M\subseteq \mathcal{E}_{1,1+\eta}(M',1). 
\end{equation}
Since $ \eta $ is arbitrary, Remark \ref{inclucomp} implies that \eqref{inclusob} holds with compact canonical injections. 
We also see that, for any given integer $ i\geq 0 $, the map $ u\mapsto u^{(i)}(0) $ is a continuous linear form on the Hilbert space $ \mathcal{H}_M $. Thus, there exists an element $ e_i $ of $ \mathcal{H}_M $ such that 
\begin{equation*}
u^{(i)}(0)=\langle e_i\vert u\rangle\ \textrm{ for any } u\in \mathcal{H}_M.
\end{equation*}

Given an element $ g $ of $ \mathcal{H}_M $, we now study the elements $ u $ of $ \mathcal{H}_M $ solving the system
\begin{equation}\label{sys1}
u^{(i)}(0)=g^{(i)}(0)\ \textrm{ for } 0\leq i<k.
\end{equation}
Denote by $\mathcal{V}_k $ the subspace of $ \mathcal{H}_M $ spanned by
$ e_0,e_1,\ldots, e_{k-1} $. If $ u $ is an element $ \sum_{j=0}^{k-1}\xi_{j}e_j $ of $ \mathcal{V}_k $, this system can be rewritten as a $ k\times k $ linear system with respect to the unknowns $ \xi_0,\ldots,\xi_{k-1} $, namely
\begin{equation*}
\sum_{j=0}^{k-1}\langle e_i\vert e_j\rangle\xi_j=g^{(i)}(0)\  \textrm{ for } 0\leq i<k. 
\end{equation*}
Taking scalar products with monomials, it is easy to see that the $ e_j $'s are linearly independent in $ \mathcal{H}_M $. Thus, $ \dim\mathcal{V}_k=k $ and the Gram matrix $ \big(\langle e_i\vert e_j\rangle\big)_{0\leq i,j<k} $ is invertible. We derive that \eqref{sys1} has a solution $ g_k $ which can be written $ \sum_{j=0}^{k-1}\xi_{j,k}e_j $ for some suitable 
$ \xi_{j,k} $ depending linearly on the derivatives $ g^{(i)}(0) $ with $ 0\leq i<k $. Rearranging, we obtain a family $ (u_{j,k})_{0\leq j<k} $ of elements of $ \mathcal{H}_M $, not depending on $ g $, such that 
\begin{equation}\label{xpres1}
g_k(t)=\sum_{j=0}^{k-1}u_{j,k}(t)g^{(j)}(0)\ \textrm{ for any } t\in I.
\end{equation}
Remark that any other solution $ u $ of \eqref{sys1} satisfies $ \langle e_i\vert g_k-u \rangle=g_k^{(i)}(0)-u^{(i)}(0)=0 $ for $ 0\leq i<k $. Thus, we have $ g_k-u \in \mathcal{V}_k^\bot $, whereas $ g_k \in \mathcal{V}_k $. Pythagoras's theorem then shows 
that $ g_k $ is the minimal solution\footnote{The variational argument in Carleman's proof is precisely used to construct a solution of $ \eqref{sys1} $ minimizing a functional $ u\mapsto I_k(u) $ which, roughly speaking, satisfies $ \lim_{k\to\infty}I_k(u)=\Vert u\Vert^2 $ in our notation.} of \eqref{sys1} in $ \mathcal{H}_M $. In particular, we have $\Vert g_k\Vert\leq \Vert g\Vert $ for any $ k\geq 1 $, and the sequence $ (g_k)_{k\geq 1} $ is bounded in $ \mathcal{H}_M $. We claim that
\begin{equation}\label{kkt}
(g_k)_{k\geq 1} \textrm{ converges to } g \textrm{ in } \mathcal{E}_{1,1+\eta}(M',1). 
\end{equation}
Indeed, since the inclusion $ \mathcal{H}_M\hookrightarrow \mathcal{E}_{1,1+\eta}(M',1) $ is compact, it suffices to check that $ g $ is the only possible limit for any subsequence of $ (g_k)_{k\geq 1} $ converging in $\mathcal{E}_{1,1+\eta}(M',1) $. This can be proved as follows: let $ h $ be the limit of such a subsequence. Taking limits in \eqref{sys1}, we immediately get 
\begin{equation}\label{kkr}
h^{(i)}(0)=g^{(i)}(0)\ \textrm{ for any } i\geq 0.
\end{equation}
The quasianalyticity assumption for $ \mathcal{E}_1(M) $ immediately implies that the local ring $ \mathcal{E}_1(M') $ is also quasianalytic, as well as its translates at any point of $ I $. Thus, \eqref{kkr} yields $ h=g $, as desired.

Now, let $ f $ be an element of $ \mathcal{E}_1(M) $. For any real $ x $ sufficiently close to $ 0 $, the function $ f_x $ defined by $ f_x(t)=f(xt) $ belongs to $ \mathcal{E}_{1,1-\eta}(M,1) $, hence to $ \mathcal{H}_M $ by \eqref{inclusob}. Applying \eqref{xpres1} and \eqref{kkt} to $ g=f_x $, and putting $ t=1 $ and $ \omega_{j,k}=j!u_{j,k}(1) $, we finally obtain the representation formula
\begin{equation}\label{repres}
f(x)=\lim_{k\to\infty}\sum_{j=0}^{k-1}\omega_{j,k}\frac{f^{(j)}(0)}{j!}x^j
\end{equation}
for $ x $ sufficiently close to $ 0 $.

\begin{rem}
Refining slightly the compactness argument above, it is possible to show that the convergence in \eqref{repres} is uniform with respect to $ x $ in a sufficiently small neighborhood of $ 0 $.
\end{rem}

\subsubsection*{Second part}
Applying \eqref{repres} to the monomial $ f(x)=x^j $ for any given $ j\geq 0 $, we get $ \lim_{k\to \infty}\omega_{j,k}=1 $. It is therefore possible to select, by induction, an increasing sequence $ (k_p)_{p\geq 0} $ of positive integers such that
\begin{equation}\label{kk1}
\sum_{j=0}^{k_{p-1}}\left\vert\omega_{j,k_p}-1\right\vert M_j\leq 1\ \textrm{ for any }p\geq 1. 
\end{equation}
Assume $ \mathcal{O}_1\subsetneq \mathcal{E}_1(M) $. Remark \ref{nonan} then implies
\begin{equation}\label{kk2}
\lim_{p\to\infty}(M_{k_p})^{1/k_p}=\infty.
\end{equation}
Consider the element $ F=\sum_{j\geq 0}F_jx^j $ of $ \mathcal{F}_1(M) $ defined by 
$ F_j=M_j $ for $ j\in\{k_p : p\geq 0\} $ and $ F_j=0 $ otherwise. We claim that $ F $ does not belong to the range of the Borel map $ T_0 $. Indeed, assume that there exists an element $ f $ of $ \mathcal{E}_1(M) $ such that $ T_0f=F $. Being given a sufficiently small number $ a>0 $, the representation formula \eqref{repres} then yields 
\begin{equation*}
f(a)=\lim_{p\to\infty}\sum_{j=0}^{k_p-1}\omega_{j,k_p}F_ja^j.
\end{equation*}
By the definition of $ F $, we have $ \sum_{j=0}^{k_p-1}\omega_{j,k_p}F_ja^j=\sum_{j=0}^{k_{p-1}}\omega_{j,k_p}F_ja^j $ and $ \sum_{j=0}^{k_{p-1}}F_ja^j=\sum_{q=0}^{p-1}M_{k_q}a^{k_q} $. Thus, we get 
\begin{equation}\label{kk3}
f(a)=\lim_{p\to\infty}\bigg(\sum_{q=0}^{p-1}M_{k_q}a^{k_q}+\sum_{j=0}^{k_{p-1}}(\omega_{j,k_p}-1)F_ja^j\bigg).
\end{equation}
By \eqref{kk1}, the second sum in the right-hand side of \eqref{kk3} is bounded by $ 1 $ uniformly with respect to $ p $. Thus, \eqref{kk3} implies that the sequence of partial sums of the series $ \sum_{q\geq 0} M_{k_q}a^{k_q} $ is bounded, which obviously contradicts \eqref{kk2}. The proof is complete.

\begin{rem}
As pointed out by Carleman, the above proof yields, in fact, a more precise statement: \textit{with any quasianalytic local ring $\mathcal{E}_n(M) $ one can associate an increasing sequence $ (k_p)_{p\geq 0} $ of positive integers such that the only lacunary series $ \sum_{p\geq 0} F_{k_p}x^{k_p} $ belonging to $ T_0\mathcal{E}_1(M) $ are the convergent ones}. One can therefrom see that $ \mathcal{F}_1(M) $ contains elements that do not even belong to $ T_0\mathcal{E}_1(N) $ for some larger quasianalytic ring $ \mathcal{E}_1(N) $.
\end{rem}

\section{The Weierstrass division property}

\subsection{Basic facts}
Let $ \mathcal{R}_n $ be one the local rings $ \mathcal{O}_n $, $ \mathcal{E}_n $, $ \mathcal{E}_n(M) $, $ \mathcal{F}_n $, $ \mathcal{F}_n(M) $. For $ x=(x_1,\ldots,x_n) $, write $ x'=(x_1,\ldots,x_{n-1}) $, so that $ x=(x',x_n) $. An element $ \varphi $ of 
$ \mathcal{R}_n $ is \emph{regular of order $ d $ with respect to $ x_n $} if $ \varphi(0,x_n)=u(x_n)x_n^d $ where $ u $ is a unit in $ \mathcal{R}_1 $. Being given such a $ \varphi $, we say that \emph{Weierstrass division holds in $\mathcal{R}_n $ for the divisor $ \varphi $} if, for any element $ f $ of $ \mathcal{R}_n $, one can find $ g $ in $ \mathcal{R}_n $ and $ h_0,\ldots,h_{d-1} $ in $ \mathcal{R}_{n-1} $ such that 
\begin{equation*}
f=\varphi g+h\quad\textrm{with}\quad h(x',x_n)=\sum_{j=0}^{d-1}h_j(x')x_n^j.
\end{equation*}

It is a classical fact\footnote{Interestingly enough, it seems that the property does not appear explicitly in the work of Weierstrass: see, in \cite{Grauert}, the historical account of contributions by Stickelberger, Sp\"ath, R\"uckert.} that $ \mathcal{O}_n $ and $ \mathcal{F}_n $ have the Weierstrass division property \cite{Grauert, Malgrange, Tougeron}. The same statement for $ \mathcal{E}_n $ is a famous result of Malgrange \cite{Malgrange}, based on a delicate study of analytic sets. Quite different proofs were given later by \L ojasiewicz, Mather, Nirenberg: see \cite{Tougeron} and the references therein. 

The most basic fact in the case of $ \mathcal{E}_n(M) $ and $ \mathcal{F}_n(M) $ can be stated as follows. Note that it does not depend, in any way, on quasianalyticity properties.

\begin{prop}\label{nodiv} 
Assume $ \mathcal{O}_n\subsetneq \mathcal{E}_n(M) $. Then the rings $ \mathcal{E}_n(M) $ and $ \mathcal{F}_n(M) $ do not have the Weierstrass division property.
\end{prop} 

\begin{proof}
We can assume $ n=2 $, $ (x_1,x_2)=(x,y) $. Consider the case of germs first. Put $ \varphi(x,y)=y^2+x $. Then $ \varphi$ is regular of order $ 2 $ with respect to $ y $. Put $ f(x,y)=\theta(y)$, where $ \theta $ comes from the statement of Theorem \ref{cartan}. We claim that in any division identity 
$f(x,y)=\varphi(x,y)g(x,y)+yh_1(x)+h_0(x)$ in $ \mathcal{E}_2 $, the germ $ h_0 $ cannot belong to $ \mathcal{E}_1(M) $. Indeed, putting $ x=-y^2 $, we get $ \theta(y)=yh_1(-y^2)+h_0(-y^2) $. Thus, $ h_0(-y^2) $ is the even part of $ \theta(y) $. We derive 
\begin{equation*}
T_0h_0(x)=\sum_{j\geq 0}(-1)^j\frac{\theta^{(2j)}(0)}{(2j)!}x^j. 
\end{equation*}
We know that $ \vert\theta^{(2j)}(0)/(2j)!\vert\geq M_{2j} $.
The logarithmic convexity of $ M $ implies $ M_{2j}\geq M_j^2 $ and the assumption $ \mathcal{O}_n\subsetneq \mathcal{E}_n(M) $ yields $ \sup_{j\geq 1}(M_j)^{1/j}=\infty $ by Corollary \ref{nonana}. We have therefore $ \sup_{j\geq 1}(M_{2j}/M_j)^{1/j}=\infty $. This shows that $ T_0h_0 $ does not belong to $ \mathcal{F}_1(M) $, hence $ h_0 $ does not belong to $ \mathcal{E}_1(M) $. The same argument also works in the formal case.
\end{proof}

\begin{rem}\label{lossofreg}
The preceding proof relies on a loss of smoothness in the division process, in the sense that the result $ h_0 $ lies in a strictly larger class of germs than the data $ f $. In sufficiently regular non-quasianalytic situations, this loss of smoothness can be controlled by the order; see \cite{ChtCht0}. Even more precise estimates exist in the formal case; see \cite{Mouze}, Annex B, Theorem 2. However, a better formal result is obtained when the usual assumption of regularity of order $ d $ with respect to $ x_n $ is replaced by a stronger one, as we shall see now.
\end{rem}

\subsection{A positive result in the formal case}
Say that a formal power series $ F=\sum_{J\in\mathbb{N}^n} F_Jx^J $ is \emph{strictly regular of order $ d $} (or \emph{$d$-regular}) with respect to $x_n$ if it satisfies
$ F_{(0,\ldots,0,d)}\neq 0 $ and $ F_J=0 $ for $ j<d $, where $ j $ denotes the length of $ J $. 

\begin{exam}
Put $n=2 $, $ (x_1,x_2)=(x,y) $. Then $ y^2+x $ is not strictly regular with respect to $ y $, whereas $ y^2+x^2 $ is strictly regular of order $ 2 $. 
\end{exam}

\begin{rem}
It is not difficult to show that any non-zero formal power series can be made strictly regular after a linear change of variables.
\end{rem}

The following theorem is due to Chaumat and Chollet \cite{ChtCht1}.

\begin{thm}\label{formalnoether}
The following properties are equivalent:\\
(i) The ring $ \mathcal{F}_n(M) $ is stable under derivation,\\
(ii) Weierstrass division holds in $ \mathcal{F}_n(M) $ for strictly regular divisors,\\
(iii) The ring $ \mathcal{F}_n(M) $ is noetherian.
\end{thm}

The implication $(ii)\Longrightarrow (iii)$ follows the inductive scheme which is classically used to derive the noetherianity of $ \mathcal{O}_n $ or $ \mathcal{F}_n $ from the Weierstrass division theorem: see, for instance, \cite{Grauert}. The difficulty of the result lies in the other implications. Remark, however, that a much simpler proof of $(i)\Longrightarrow (ii)$ can be found in \cite{Mouze}, Annex B, Theorem 1. 

\subsection{The case of non-quasianalytic germs}\label{strong}
As we have seen in Proposition \ref{nodiv}, the ring $ \mathcal{E}_n(M) $ does not have the general Weierstrass division property. It is thus natural to ask whether division is possible for ``good''  divisors, as this happens with $ \mathcal{F}_n(M) $. In non-quasianalytic situations, there are actually some (rather particular) results. 

Say that the logarithmically convex sequence $ M $ is \emph{strongly regular} if it satisfies the strong non-quasianalyticity condition \eqref{SNQA} together with the \emph{moderate growth condition}
\begin{equation}\label{modgr}
\sup_{j\geq 1,\, k\geq 1}\left(\frac{M_{j+k}}{M_jM_k}\right)^{\frac{1}{j+k}}<\infty.
\end{equation}
Condition \eqref{modgr} is obviously stronger than the condition of stability under derivation given in Corollary \ref{stabder}. It has an interpretation in terms of stability under the action of so-called ultradifferential operators \cite{Komatsu}. Let us also mention that strong regularity ensures that a version of Whitney's extension theorem holds for the corresponding classes of ultradifferentiable functions and jets: see, for instance, \cite{BBMT}.
A typical example of strongly regular sequence is given by the Gevrey sequences $ M_j=(j!)^\alpha $ with $ \alpha>0 $. The $ q $-Gevrey sequence $ M_j =q^{j^2} $ with $ q>1 $ does not satisfy the moderate growth condition, hence is not strongly regular (the corresponding ring $ \mathcal{E}_n(M) $ is merely stable under derivation and strongly non-quasianalytic).

As suggested in Remark \ref{lossofreg}, in strongly regular classes, Weierstrass division holds with a loss of smoothness related to the order \cite{ChtCht0}. Here, we shall rather present a particular situation where division is possible without loss of smoothness. 

We need first to recall a few facts about regular separation. Let $ \varphi $ be a real analytic function germ. Denote by $ Z_\varphi $ the zero set of the natural complexification of $ \varphi $. This is a germ of complex analytic set. Denote by $ X_\varphi $ the real zero set of $ \varphi $, that is, 
$ Z_\varphi\cap \mathbb{R}^n $ (we always view $ \mathbb{R}^n $ as a totally real subset of $ \mathbb{C}^n $ in the natural way). Denote by $ \mathcal{T} $ the set of real numbers $ \tau $ for which one can find a constant $ c>0 $ such that the inequality $ \dist(x,Z_\varphi)\geq c\dist(x,X_\varphi)^\tau $ holds for any $ x $ in a neighborhood of $ 0 $ in $ \mathbb{R}^n $. Then $ \mathcal{T} $ is non-void \cite{Bochnak-Risler}, and bounded below by $ 1 $. The number $ \tau(\varphi)=\inf\mathcal{T} $ is called the \emph{\L ojasiewicz exponent for the regular separation between $ Z_\varphi $ and $ \mathbb{R}^n $}. As shown in \cite{Bochnak-Risler}, the exponent $ \tau(\varphi) $ is rational and it belongs to $ \mathcal{T} $. 

\begin{prop}\label{simplediv}
Assume that $ M $ is strongly regular and let $ \varphi $ be a real-analytic function germ such that $ X_\varphi=\{0\} $ and $ \tau(\varphi)=1 $. 
Then Weierstrass division holds in $\mathcal{E}_n(M) $ for the divisor $ \varphi $.
\end{prop}

This proposition can be readily derived from Theorems \ref{T0surj}, \ref{formalnoether} and \ref{ellipt} of the present article. Geometrically speaking, the assumptions mean that the divisor has an isolated real zero at the origin and that its complex zero set lies outside a conical neighborhood of $ \mathbb{R}^n\setminus\{0\} $ in $ \mathbb{C}^n $.

\begin{exam}\label{square}
Put $ n=2 $, $ (x_1,x_2)=(x,y) $. According to Proposition \ref{simplediv}, if the sequence $ M $ is strongly regular, Weierstrass division holds in $\mathcal{E}_2(M) $ for the divisor $ \varphi(x,y)=y^2+x^2 $. Quite surprisingly, this is no longer true in the quasianalytic case, as we shall see now.
\end{exam}

\subsection{The quasianalytic case}\label{hyperb}
In this case, there is a fairly complete solution of the division problem in the typical situation of distinguished polynomials in $ \mathcal{E}_{n-1}(M)[x_n] $, that is, for divisors $ \varphi $ of the form 
\begin{equation}\label{poly}
\varphi(x)=x_n^d+a_1(x')x_n^{d-1}+\cdots+a_d(x') 
\end{equation}
with $a_j\in \mathcal{E}_{n-1}(M)$ and $ a_j(0)=0 $ for $ 1\leq j\leq d $. Such a polynomial  $ \varphi $ is said to be \emph{hyperbolic} if there exists a neighborhood $ U $ of $ 0 $ in $ \mathbb{R}^{n-1} $ such that for any $ x' \in U $, all the roots of $ \varphi(x',\cdot) $ are real. The following theorem is due to Childress \cite{Childress}. 

\begin{thm}\label{child}
Assume that the local ring $ \mathcal{E}_n(M) $ is quasianalytic, with $ \mathcal{O}_n\subsetneq\mathcal{E}_n(M) $, and that Weierstrass division holds in $\mathcal{E}_n(M) $ for a distinguished polynomial $ \varphi $. Then $ \varphi $ is hyperbolic. 
\end{thm}

In particular, as announced in Example \ref{square}, we see that in the quasi-analytic case, Weierstrass division does not hold in $ \mathcal{E}_2(M) $ for the divisor $ \varphi(x,y)=y^2+x^2 $. On the positive side, we have the following result, due to Chaumat and Chollet \cite{ChtCht}. 

\begin{thm}\label{chch}
Assume that the local ring $ \mathcal{E}_n(M) $ (either quasianalytic or not) is stable under derivation. Let $ \varphi $ be a hyperbolic polynomial as above. Then Weierstrass division holds in $ \mathcal{E}_n(M) $ for the divisor $ \varphi $. 
\end{thm}

Gathering Theorems \ref{child} and \ref{chch}, we get the following necessary and sufficient condition.

\begin{cor}
Assume that $ \mathcal{E}_n(M) $ is quasianalytic, stable under derivation, and $ \mathcal{O}_n\subsetneq \mathcal{E}_n(M) $. Let $ \varphi $ be a distinguished polynomial in $\mathcal{E}_{n-1}(M)[x_n] $. Then Weierstrass division holds in $ \mathcal{E}_n(M) $ for the divisor $ \varphi $ if and only if $ \varphi $ is hyperbolic. 
\end{cor}

\subsection{A proof of Theorem \ref{child}}\label{prereq}
We shall present a complete proof of Theorem \ref{child}. We do not pretend originality: all the arguments come from \cite{Childress} and \cite{ChtCht}. However, we think they are worth being reproduced here, since they are quite typical of our topic. This is also an opportunity to introduce additional material that will be used later. 

\subsubsection*{Topological prerequisites}
We shall need the Banach spaces $ \mathcal{E}_{n,\sigma}(M,\nu) $ defined in \S \ref{revisit}. For $ \nu'>\nu $ and $ \sigma'>\sigma $, we have a commutative diagram
\begin{equation*}
\xymatrix{
\mathcal{E}_{n,\sigma}(M,\nu) \ar[d] \ar[r] & \mathcal{E}_{n,\sigma}(M,\nu') \ar[d]\\
\mathcal{E}_{n,\sigma'}(M,\nu) \ar[r] &\mathcal{E}_{n,\sigma'}(M,\nu')
}
\end{equation*}
where the horizontal arrows are given by the natural restriction maps, and the vertical arrows by canonical injections. The quasianalyticity assumption implies that the restriction maps are injective. The canonical injections are compact maps, as mentioned  in Remark \ref{inclucomp}. Put $ \sigma =\nu $, $ \sigma'=\nu'=\nu+1 $, $ E_\nu=\mathcal{E}_{n,\nu}(M,\nu) $, $ \Vert\cdot\Vert_\nu=\Vert\cdot\Vert_{\nu,\nu} $ and notice that we have $ \mathcal{E}_n(M)=\bigcup_{\nu\geq 1}E_\nu $ as a set. From the considerations above, we get a compact injection $ E_\nu\longrightarrow E_{\nu+1} $. Thus, $ \mathcal{E}_n(M) $ is a so-called \emph{Silva space}, or \emph{(DFS)-space}, for the corresponding inductive topology.  

We shall not enumerate the properties of Silva spaces. The interested reader will find more information in a very readable form in \cite{Grauert}, Chapter I, \S 7--8. One can also refer to \cite{Meise-Vogt}, Chapters 24 and 25.  

\subsubsection*{Scheme of the proof}
It proceeds by \textit{reductio ad absurdum}. We assume that\\
(i) Weierstrass division holds in $ \mathcal{E}_n(M) $ for the divisor $ \varphi $,\\
(ii) the polynomial $ \varphi $ is not hyperbolic,\\
and we look for a contradiction in two steps. Assumption (i) will be used in the first step, assumption (ii) in the second. 

\subsubsection*{First step}
Consider the map $ A : \mathcal{E}_n(M)\times\left(\mathcal{E}_{n-1}(M)\right)^d\longrightarrow \mathcal{E}_n(M) $ defined by
\begin{equation*} 
A(g,h_0,\ldots,h_{d-1})=\varphi g+h \ \textrm{ with }\   h(x',x_n)=\sum_{j=0}^{d-1}h_j(x')x_n^j .
\end{equation*} 
It is not difficult to check that $ A $ is a linear continuous operator between the product of Silva spaces $ \mathcal{E}_n(M)\times\left(\mathcal{E}_{n-1}(M)\right)^d $ and the Silva space $ \mathcal{E}_n(M) $. Moreover, it is known that the quotient and the remainder of Weierstrass division identities in $ \mathcal{F}_n $ are unique \cite{Malgrange, Tougeron}. This, together with the quasianalyticity of $ \mathcal{E}_n(M) $, implies that $ A $ is injective. It is also surjective by assumption (i) above. Thus, by the open mapping theorem in the Grothendieck version (\cite{Meise-Vogt}, Theorem 24.30 and Remark 24.36), the map $ A $ is a topological isomorphism. The continuity of $ A^{-1} $ can be written as follows: for any $ \nu \geq 1 $, there exist an integer $ \lambda_\nu\geq 1 $ and a real constant $ C_\nu>0 $ such that 
\begin{equation*}
\Vert A^{-1}f\Vert_{\lambda_\nu}\leq C_\nu \Vert f\Vert_\nu \ \textrm{ for any } f\in E_\nu, 
\end{equation*}
which means that any $ f $ in $ E_\nu $ can be written $ f(x)=\varphi(x) g(x) +\sum_{j=0}^{d-1}h_j(x')x_n^j $ for $ x=(x',x_n)\in I_{\lambda_\nu}^n $, with \textit{a priori} estimates
\begin{equation}\label{apriori}
\Vert g\Vert_{\lambda_\nu}\leq C_\nu\Vert f\Vert_\nu\ \textrm{ and }\ \Vert h_j\Vert_{\lambda_\nu}\leq C_\nu\Vert f\Vert_\nu\ \textrm{ for } 0\leq j<d.
\end{equation}  

\subsubsection*{Second step}
We now use an argument due to Chaumat and Chollet \cite{ChtCht}, which is simpler than the original proof of Childress. Consider any one-variable polynomial $ P $ in $\mathbb{C}[z] $, and write the euclidean division of $ P $ by the generic polynomial of degree $ d $, that is, by $ F_d(\mu,z)=z^d+\mu_1 z^{d-1}+\cdots+\mu_d $ with $ \mu=(\mu_1,\ldots,\mu_d)\in \mathbb{C}^d $. We have
\begin{equation}\label{basicdiv}
P(z)=F_d(\mu,z)G(\mu,z)+H(\mu,z) \ \textrm{ with }\ H(\mu,z)=\sum_{j=0}^{d-1}H_j(\mu)z^j.
\end{equation}
Both terms $ G(\mu,z) $ and $ H(\mu,z) $ are polynomials in $ \mathbb{C}[z] $ with coefficients depending on $ P $ and $ \mu $. Now, for any $ x' $ in a sufficiently small neighborhood of $ 0 $ in $ \mathbb{R}^{n-1} $, put $ \mu=a(x')=(a_1(x'),\ldots,a_d(x')) $, where the $ a_j $'s are the coefficients of $ \varphi $ in \eqref{poly}. Put also $ \widetilde{G}(x',z)=G(a(x'),z) $ and $ \widetilde{H}_j(x',z)=H_j(a(x'),z)$ for $ 0\leq j<d $. The division identity \eqref{basicdiv} becomes 
\begin{equation}\label{diveuc}
P(z)=\varphi(x',z)\widetilde{G}(x',z)+\sum_{j=0}^{d-1}\widetilde{H}_j(x')z^j.
\end{equation}
Now, write the Weierstrass division property for the function $ f $ given by $ f(x',x_n)=P(x_n) $. For $ x $ in $ I_{\lambda_\nu}^n $, we have
\begin{equation}\label{div2}
P(x_n)=\varphi(x',x_n)g(x',x_n)+\sum_{j=0}^{d-1}h_j(x')x_n^j.
\end{equation}
From \eqref{diveuc}, \eqref{div2} and the uniqueness of euclidean division in $ \mathbb{C}[x_n] $, we derive 
\begin{equation}\label{ident}
g=\widetilde{G} \ \textrm{ and }\ h_j=\widetilde{H_j} \ \textrm{ for }\ 0\leq j<d.
\end{equation}
Now, we use the non-hyperbolicity assumption (ii). Fix an integer $ \nu\geq 1 $. Then there exist $ x' $ in $ I_{\lambda_\nu}^{n-1} $ and $ z_0 $ in $ \mathbb{C}\setminus\mathbb{R} $ such that $ \varphi(x',z_0)=0 $. From \eqref{diveuc} and \eqref{ident}, we derive 
\begin{equation*}
P(z_0)=\sum_{j=0}^{d-1}h_j(x')z_0^j.
\end{equation*}
The \textit{a priori} estimates \eqref{apriori} then yield
\begin{equation}\label{atz}
\vert P(z_0)\vert\leq C_\nu\bigg(\sum_{j=0}^{d-1}\vert z_0\vert^j\bigg)\Vert f\Vert_\nu
\ \textrm{ with }\ \Vert f\Vert_\nu=\sup_{j\geq 0,\, \vert x_n\vert<1/\nu} 
\frac{\vert P^{(j)}(x_n)\vert}{\nu^jj!M_j}.
\end{equation}
For any $ x_n $ in $ I_\nu $, we majorize $ \vert P^{(j)}(x_n)\vert $ by means of the Cauchy formula used on the disc of center $ x_n $ and radius $ \vert \Im z_0\vert/2 $.  
Putting $W=\{z\in\mathbb{C} : \dist(z,I_\nu)<\vert\Im z_0\vert/2 \} $, we obtain
\begin{equation}\label{Cauchy} 
\vert P^{(j)}(x_n)\vert\leq j!(2/\vert\Im z_0\vert)^j\sup_{\overline{W}}\vert P\vert.
\end{equation}

Now, we use the assumption $ \mathcal{O}_n\subsetneq\mathcal{E}_n(M) $. Remark \ref{nonan} yields $ \lim_{j\to \infty} (M_j)^{1/j}=\infty $, hence the quantity $ D_\nu $ defined by 
\begin{equation*}
D_\nu=\sup_{j\geq 0}(2/\vert\Im z_0\vert)^j
\frac{1}{M_j} 
\end{equation*}
is finite. From \eqref{atz} and \eqref{Cauchy}, we then obtain
\begin{equation}\label{final}
\vert P(z_0)\vert\leq C\sup_{\overline{W}}\vert P\vert\ \textrm{ for any } P\in\mathbb{C}[z],
\end{equation}
with $ C=C_\nu D_\nu\sum_{j=0}^{d-1}\vert z_0\vert^j $. 
We claim that \eqref{final} is impossible. Indeed, put $ K=\overline{W}\cup\{z_0\} $. Then $ K $ is a compact subset of $ \mathbb{C}$ and $ \mathbb{C}\setminus K $ is connected. Thanks to Runge's approximation theorem (\cite{Remmert}, Chapter 12), we can construct a sequence of polynomials $ (P_k)_{k\geq 1} $ such that $ \lim_{k\to\infty} P_k(z_0)=1 $ and $ \lim_{k\to\infty} \sup_{\overline{W}}\vert P_k\vert=0 $. Obviously, this contradicts \eqref{final}. The proof is complete. 

\subsection{An open problem}\label{open}
It is well-known that the rings $ \mathcal{O}_n $ and $ \mathcal{F}_n $ are noetherian: this is easy to see for $ n=1 $, and the proof for higher dimensions goes by induction, \textit{via} the Weierstrass division theorem \cite{Grauert}. It is also clear that $ \mathcal{E}_n $ is not noetherian: indeed, if $ \imax $ denotes its maximal ideal, we have $ \bigcap_{k\geq 1}\imax^k\neq \{0\} $, since this intersection contains all flat germs\footnote{A smooth function germ is said to be \emph{flat} if it vanishes at $ 0 $, together with all its derivatives.}. Thus, the ring $ \mathcal{E}_n $ does not satisfy the conclusion of Krull's intersection theorem (for which we refer to \cite{Grauert}, Appendix, \S 2). The same argument shows that in the non-quasianalytic case, the ring $ \mathcal{E}_n(M) $ is not noetherian. On the positive side, we have learned from Theorem \ref{formalnoether} that $ \mathcal{F}_n(M) $ is noetherian if and only if it is stable under derivation. We are thus led to the following natural question: 

\begin{prob}\label{tobeornot}
Assume that the ring $ \mathcal{E}_n(M) $ is quasianalytic, stable under derivation, and that $ \mathcal{O}_n\subsetneq \mathcal{E}_n(M) $. Is it then a noetherian ring?
\end{prob}

The failure of Weierstrass division described in Theorem \ref{child} shows that it is not possible to mimic the classical proof of the noetherianity of $ \mathcal{O}_n $ or $ \mathcal{F}_n $. However, this is not enough to say that noetherianity fails.   
Problem \ref{tobeornot} is, as of this writing, still open. It is obviously an important question, since noetherianity has a number of algebraic and geometric consequences. In what follows, we suggest two possible, albeit so far inconclusive, directions to study the problem.  

\subsubsection*{First direction: flatness}
In the quasianalytic case, the Borel map $ T_0 $ embeds $ \mathcal{E}_n(M) $ as a subring of the ring $ \mathcal{F}_n $ of all formal power series, and an element of $ T_0\mathcal{E}_n(M) $ is invertible in $ T_0\mathcal{E}_n(M) $ if and only if it is invertible in $ \mathcal{F}_n $. The noetherianity of such a subring of $ \mathcal{F}_n $ can be reduced to the study of a flatness property of modules\footnote{A module $ M $ over a ring $ R $ is said to be \emph{flat} if, for any ideal $ \mathcal{I} $ of $ R $, the natural map $ \mathcal{I}\otimes_R M \rightarrow M $ is injective. We refer the reader to \cite{Malgrange} for more details and equivalent formulations.} by means of the following lemma. 

\begin{lem}
Let $ \mathcal{C}_n $ be a local subring of $\mathcal{F}_n $ containing  $ \mathcal{O}_n $. Assume that the maximal ideal $\imax $ of $ \mathcal{C}_n $ coincides with $ \hat\imax\cap\mathcal{C}_n $, where $ \hat\imax $ denotes the maximal ideal of $ \mathcal{F}_n $. Then $ \mathcal{C}_n $ is noetherian if and only if $ \mathcal{F}_n $ is a flat $ \mathcal{C}_n $-module.
\end{lem}

\begin{proof}
The proof of sufficiency is elementary. First, a basic exercise yields the following variation on the usual definition of noetherianity: $ \mathcal{C}_n $ is noetherian if and only if any increasing sequence of \emph{finitely generated} ideals $ \mathcal{I}_0\subseteq \mathcal{I}_1\subseteq\cdots\subseteq \mathcal{I}_j\subseteq \mathcal{I}_{j+1}\subseteq \cdots $ of $ \mathcal{C}_n $ stabilizes. Now, assume that $ \mathcal{F}_n $ is flat over $ \mathcal{C}_n $ and consider such a sequence. Denote by $ a_1^{(j)},\ldots,a_{p_j}^{(j)} $ a system of generators of $ \mathcal{I}_j $. The sequence of ideals $ \mathcal{I}_j\mathcal{F}_n $ increases in the noetherian ring $ \mathcal{F}_n $, hence there exists an integer $ j_0\geq 0 $ such that $ \mathcal{I}_j\mathcal{F}_n=\mathcal{I}_{j_0}\mathcal{F}_n $ for any $ j\geq j_0 $. In particular, for any $ j\geq j_0 $ and any $ k=1,\ldots,p_j $, we have 
$ a^{(j)}_k=\sum_{i=1}^{p_{j_0}}f^{(j)}_{ki}a^{(j_0)}_i $ for some suitable elements $ f^{(j)}_{ki} $ of $\mathcal{F}_n $. In other words, the $ (p_{j_0}+1) $-tuple $ r=(1,-f^{(j)}_{k1},\ldots,-f^{(j)}_{kp_{j_0}}) $ provides a linear relation between $ a^{(j)}_k,a^{(j_0)}_1,\ldots,a^{(j_0)}_{p_{j_0}} $ with coefficients in $ \mathcal{F}_n $.  The flatness assumption then implies the existence of an identity $ r=\sum_{q=1}^m g_q r_q $ where each $ g_q $ belongs to $ \mathcal{F}_n $ and each $ r_q $ provides a linear relation between $ a^{(j)}_k,a^{(j_0)}_1,\ldots,a^{(j_0)}_{p_{j_0}} $ with coefficients $ r_{qi}$ in $ \mathcal{C}_n $. The identity implies $ 1=\sum_{q=1}^m g_q r_{q1} $, so that, for at least one index $ q $, the coefficient $ r_{q1} $ must be invertible in $ \mathcal{F}_n $, hence in $ \mathcal{C}_n $ by the assumption on maximal ideals. The relation 
$  r_{q1}a^{(j)}_k+\sum_{i=2}^{p_{j_0}+1}r_{qi}a^{(j_0)}_k=0 $
then shows that $ a^{(j)}_k $ belongs to $ \mathcal{I}_{j_0} $. This yields $ \mathcal{I}_j=\mathcal{I}_{j_0} $ for $ j\geq j_0 $, hence the noetherianity. 

Necessity can be obtained as follows. Since $ \imax=\hat\imax\cap\mathcal{C}_n $, the $\imax$-adic topology of $ \mathcal{C}_n $ is induced by the $\hat\imax $-adic topology of $ \mathcal{F}_n $. Therefore, considering the inclusion $ \mathcal{O}_n \subseteq \mathcal{C}_n \subseteq\mathcal{F}_n $ and taking completions, we get $ \widehat{\mathcal{C}_n}=\mathcal{F}_n $. The conclusion follows since the completion of a noetherian ring is flat over it: see, for instance, \cite{Malgrange}, Chapter III.
\end{proof}

The preceding lemma insists on the algebraic flavor of problem \ref{tobeornot}. We shall see hereafter how to relate it to considerations of analysis.  

\subsubsection*{Second direction: closedness}
The following result is folklore. We include a proof for the reader's convenience. 

\begin{lem}\label{closed1}
Let $ \mathcal{C}_n $ be a local subring of $\mathcal{F}_n $ whose maximal ideal $ \imax $ coincides with $ \hat\imax\cap\mathcal{C}_n $, where $ \hat\imax $ denotes the maximal ideal of $ \mathcal{F}_n $. Assume that $ \mathcal{C}_n $ is also a topological vector space and that the inclusion $ \mathcal{C}_n\hookrightarrow\mathcal{F}_n $ is continuous for the product topology\footnote{That is, the topology of simple convergence of coefficients.} on $ \mathcal{F}_n $.
If the ring $\mathcal{C}_n $ is noetherian, then all its ideals are closed.
\end{lem}

\begin{proof}
For any integer $ k\geq 0 $, denote by $ \mathcal{P}_{n,k} $ the vector space of polynomials of degree at most $ k $ in $ x_1,\ldots,x_n $. Consider the truncation operator 
$ R_k : \mathcal{C}_n\longrightarrow \mathcal{P}_{n,k} $ which, to any element of $ \mathcal{C}_n $, associates the sum of its monomials of degree at most $ k$. This map is continuous for the topology induced on $ \mathcal{C}_n $ by the product topology of $ \mathcal{F}_n $, hence for the (finer) topology of $ \mathcal{C}_n $. Now let $ \mathcal{I} $ be an ideal of $ \mathcal{C}_n $. Then $ R_k\mathcal{I} $ is a vector subspace of the finite-dimensional space $ \mathcal{P}_{n,k} $, hence closed therein. The continuity of $ R_k $ implies that $ R_k^{-1}(R_k\mathcal{I}) $ is closed in $ \mathcal{C}_n $. Now, remark that the assumption $ \imax=\hat\imax\cap\mathcal{C}_n $ implies
$ R_k^{-1}(R_k\mathcal{I})=\mathcal{I}+\imax^{k+1} $. If $\mathcal{C}_n $ is noetherian, we also have $ \mathcal{I}=\bigcap_{k\geq 0}(\mathcal{I}+\imax^{k+1}) $ by Krull's intersection theorem, and $ \mathcal{I} $ is therefore closed. 
\end{proof}

Now, consider a quasianalytic local ring $ \mathcal{E}_n(M) $, endowed with the Silva space topology described in \S \ref{prereq}. Each linear map $ f\longmapsto D^Jf(0) $ is continuous on $ \mathcal{E}_n(M) $, since we have 
$ \vert D^Jf(0)\vert\leq \Vert f\Vert_\nu\nu^jj!M_j $ for any integer $ \nu\geq 1 $ and any element $ f $ of $ E_\nu $. Thus, the Borel map identifies $ \mathcal{E}_n(M) $  with a local subring of $ \mathcal{F}_n $ satisfying the assumptions of 
Lemma \ref{closed1}. The existence of a non-closed ideal in $ \mathcal{E}_n(M) $ would then show that $ \mathcal{E}_n(M) $ is not noetherian. We do not know whether such ideals exist, although some clues make the conjecture plausible, as explained in the next section. 

\section{On closed ideals}

\subsection{A notion of closedness}
In view of \S \ref{open}, there is a strong motivation to study closed ideals in local rings $ \mathcal{E}_n(M) $. However, we have to be careful about the notion of closedness to be considered: indeed, the topological prerequisites of \S \ref{prereq} do not extend to the non-quasianalytic case, since the restriction maps $ \mathcal{E}_{n,\sigma}(M,\nu)\longrightarrow \mathcal{E}_{n,\sigma}(M,\nu') $ are no longer injective in this situation. We shall be interested in the typical case of principal ideals $ \mathcal{I}=\varphi\mathcal{E}_n(M) $ generated by a real-analytic function germ $ \varphi $. A good historical reason to study this case lies in the famous results of H\"ormander and \L ojasiewicz on the division of distributions, which, in dual form, can be stated as follows.
Let $ \varphi $ be an analytic function in a open subset $ U $ of $\mathbb{R}^n $. Then the ideal $ \varphi C^\infty(U) $ is closed for the usual Fr\'echet topology on $ C^\infty(U) $. For the details, we refer to \cite{Malgrange, Tougeron} and the references therein\footnote{Alternatively, an elementary presentation of these ideas can be found in the introductory paper \cite{Malgrange2}.}.

Now, for a given integer $ \nu\geq 1 $, consider the Banach 
spaces $ \mathcal{E}_{n,\sigma}(M,\nu) $ defined in \S \ref{revisit}. Remark \ref{inclucomp} implies that the union
$ \mathcal{E}_n(M,\nu)=\bigcup_{\sigma\in\mathbb{N}\setminus\{0\}}\mathcal{E}_{n,\sigma}(M,\nu) $, endowed with the inductive topology, is a Silva space. For $ \nu $ large enough, $ \varphi $ is analytic in $ I_\nu^n $ and we can consider the ideal $ \mathcal{I}(\nu)=\varphi\mathcal{E}_n(M,\nu) $. We shall say that $ \mathcal{I} $ is \emph{closed} if every element of the closure of $ \mathcal{I}(\nu) $ in the Silva space $ \mathcal{E}_n(M,\nu) $ defines a germ belonging to $ \mathcal{I} $. 

\begin{rem}
Consider a sequence of $ (U_\nu)_{\nu\geq 1} $ of 
open, convex, neighborhoods of $ 0 $  
such that $ \overline{U_\nu}\subseteq U_{\nu+1} $ and $ \bigcap_{\nu\geq 1}U_\nu=\{0\} $. The preceding notion of closedness does not change if the cube $ I_\nu^n $ is replaced by $ U_\nu $ in the definition of the Banach spaces $ \mathcal{E}_{n,\sigma}(M,\nu) $.
\end{rem}

\begin{rem}
Assume that $ \mathcal{E}_n(M) $ is quasianalytic. If the ideal 
$ \mathcal{I}=\varphi\mathcal{E}_n(M) $ is closed for the Silva space topology defined in \S \ref{prereq}, then it is also closed in the above sense. Indeed, the injection 
$ \iota_\nu : \mathcal{E}_n(M,\nu)\longrightarrow\mathcal{E}_n(M) $ 
is continuous, thus $ \iota_\nu^{-1}(\mathcal{I}) $ is 
closed in $ \mathcal{E}_n(M,\nu) $. Since we have obviously $ \mathcal{I}(\nu)\subseteq \iota_\nu^{-1}(\mathcal{I}) $, we derive that the closure of $ \mathcal{I}(\nu) $ in $ \mathcal{E}_n(M,\nu) $ lies in $ \iota_\nu^{-1}(\mathcal{I}) $, hence the result.
\end{rem}

\begin{rem}\label{remind}
From the preceding remark and from the considerations in \S \ref{open}, we see that the noetherianity of quasianalytic local rings $ \mathcal{E}_n(M) $ would imply the closedness of the ideal $ \varphi\mathcal{E}_n(M) $ for any given real-analytic function germ $ \varphi $. 
\end{rem}

Now, our main concern will be to describe closedness conditions for the ideal $ \mathcal{I}=\varphi\mathcal{E}_n(M) $ in terms of geometric or algebraic features of $ \varphi $. 

\subsection{The quasianalytic case}
In this case, not much is known about the closedness properties of ideals $ \mathcal{I} $ as above. We do not know an example of a non-closed ideal (as pointed out in Remark \ref{remind}, this would prove that $ \mathcal{E}_n(M) $ fails to be noetherian).
But it is also difficult to give non-trivial examples of generators $ \varphi $ for which $ \mathcal{I} $ is known to be closed.
The following statement essentially summarizes the current knowledge on the question. Notice that it is, in fact, not specific to quasianalytic situations! 

\begin{thm}\label{qaclos}
Assume that the local ring $ \mathcal{E}_n(M) $ (either quasianalytic or not) is stable under derivation. Let $ \varphi $ be a real-analytic germ at the origin of $ \mathbb{R}^n $. The ideal $ \mathcal{I} =\varphi\mathcal{E}_n(M) $ is closed in each of the following cases:\\
(i) the generator $ \varphi $ is a monomial,\\
(ii) the generator $ \varphi $ is a homogeneous polynomial with an isolated real critical point at the origin,\\
(iii) the generator $ \varphi $ is a hyperbolic polynomial.
\end{thm}

Case \textit{(i)} is simple: if $ \varphi(x)=x^J $ for some multi-index $ J $, every function $ f $ belonging to the closure of $ \mathcal{I}(\nu) $ in $ \mathcal{E}_n(M,\nu) $ can be written $ f(x)=x^Jg(x) $, where $ g(x) $ is given by the Taylor formula with integral remainder. This expression involves a derivative of $ f $; using the stability of $ \mathcal{E}_n(M) $ under derivation, it is easy to get the desired estimates. Case \textit{(ii)} is much more complicated. The result appeared in \cite{VT3} and provided the first non-trivial examples of closed principal ideals in the quasianalytic setting. Finally, case \textit{(iii)} derives from Theorem \ref{chch} in the same way as its $ C^\infty $ analogue, which can be found\footnote{We briefly recall the argument for the reader's convenience: if $ f $ belongs to the closure of $ \varphi\mathcal{E}_n(M,\nu) $ in $ \mathcal{E}_n(M,\nu) $, then we have $ f(x',x_n)=0 $ whenever $ \varphi(x',x_n)=0 $, because the evaluation maps are continuous for the topology under consideration. Moreover, by a classical argument, for any  neighborhood $ W $ of $ 0 $ in $ \mathbb{R}^{n-1}\times\mathbb{C} $, there is a neighborhood $ U$ (resp. $ V $) of $ 0 $ in $ \mathbb{R}^{n-1} $ (resp. $\mathbb{C} $) satisfying $ U\times V\subseteq W $ and such that for any $ x'\in U $, the polynomial $ \varphi(x',\cdot) $ has $ d $ roots in $ V $. The hyperbolicity assumption implies that all these roots are real. Now, by Theorem \ref{chch}, we have the division identity $ f=\varphi g +h $ where $ g $ belongs to $ \mathcal{E}_n(M) $ and $ h $ is a polynomial of degree at most $ d-1 $ in $ \mathcal{E}_{n-1}(M)[x_n] $. Choosing $ W $ small enough, we see that for any $ x' $ in $ U $, the polynomial $ h(x',\cdot) $ has $ d $ roots. We have therefore $ h=0 $, hence $ f=\varphi g $ with $ g\in\mathcal{E}_n(M) $.} in \cite{Tougeron}, Chapter V, Remark 4.7. 

In general, there seems to be no obvious link between the closedness of $ \varphi\mathcal{E}_n(M) $ and the validity of Weierstrass division by $ \varphi $, as shown by the following example. 

\begin{exam}\label{x2y4}
Put $ n=2 $, $ (x_1,x_2)=(x,y) $ and $ \varphi(x,y)=y^2+x^2 $. Assume that $ \mathcal{E}_2(M) $ is stable under derivation. Case \textit{(ii)} of Theorem \ref{qaclos} shows that the ideal $ \varphi \mathcal{E}_2(M) $ is closed. However, in the quasianalytic case, we know from \S \ref{hyperb} that Weierstrass division does not hold in $ \mathcal{E}_2(M) $ for the divisor $ \varphi $. 
\end{exam}

\begin{rem}
Although $ \varphi(x,y)=y^2+x^2 $ is a very simple polynomial, the closedness of $\varphi\mathcal{E}_2(M) $ in the preceding example is far from being trivial. One can appreciate this by comparison with $ \psi(x,y)=y^2+x^4 $: in the quasianalytic case, we do not know whether $ \psi\mathcal{E}_2(M) $ is closed. Let us now mention that, for sufficiently regular non-quasianalytic situations, $ \varphi\mathcal{E}_2(M) $ is closed, whereas $ \psi\mathcal{E}_2(M) $ is not, as we shall see in Example \ref{exzerois} below. 
\end{rem}

\subsection{The strongly regular case}\label{off}
We study now the case of strongly regular sequences, as defined in \S \ref{strong}. Although this non-quasianalytic situation is not the central theme of the present article, its specific results provide an interesting insight on the role of geometry in closedness properties. 

Whitney's spectral theorem for ultradifferentiable classes is a key step towards these results. This theorem, established in \cite{ChtCht1}, can be presented as follows.
For any integer $ \nu\geq 1 $ and any point $ a $ in $ I_\nu^n $, denote by $ T_a $ the Borel map at $ a $, that is, the map $ T_a : \mathcal{E}_n(M,\nu)\longrightarrow \mathcal{F}_n(M) $ defined by $ T_af =T_0f(\cdot +a) $. With any ideal $ \mathcal{J} $ of $ \mathcal{E}_n(M,\nu) $ we associate the ideal $ \widehat{\mathcal{J}}_a=T_a^{-1}(T_a\mathcal{J}) $ given by the functions $ f $ of $ \mathcal{E}_n(M,\nu) $ sharing their Taylor expansion at $ a $ with an element of $ \mathcal{J} $. The \textit{formal ideal} $ \widehat{\mathcal{J}} $ of $ \mathcal{J} $ is then defined by $ \widehat{\mathcal{J}}=\bigcap_{a\in I_\nu^n}\widehat{\mathcal{J}}_a $, and we have the following statement.

\begin{thm}\label{spec}
If the sequence $ M $ is strongly regular, the closure $ \overline{\mathcal{J}} $ of the ideal $ \mathcal{J} $ in the Silva space 
$ \mathcal{E}_n(M,\nu) $ 
satisfies $\overline{\mathcal{J}}=\widehat{\mathcal{J}}$.
\end{thm}

Denote by $ \imax_M^\infty $ the ideal $ \bigcap_{k\geq 1}\imax_M^k $ of flat germs in $ \mathcal{E}_n(M) $. Using Theorem \ref{spec} with $ \mathcal{J}=\mathcal{I}(\nu) $, we get an obvious corollary.

\begin{cor}\label{ellip} 
Assume that the sequence $ M $ is strongly regular and that the real-analytic germ $ \varphi $ has an isolated real zero at the origin. Then the ideal $ \mathcal{I}= \varphi\mathcal{E}_n(M) $ is closed if and only if it contains $ \imax_M^\infty $. 
\end{cor}

For any real number $ s\geq 1 $, denote by $ M^s $ the sequence $ \big((M_j)^s\big)_{j\geq 0} $. As in \S \ref{strong}, denote by $ \tau(\varphi) $ the \L ojasiewicz exponent for the regular separation between the complex zero set of $ \varphi $ and $ \mathbb{R}^n $. The following theorem comes from \cite{VT1, VT2}. 

\begin{thm}\label{ellipt}
Assume that the sequence $ M $ is strongly regular and that the real-analytic germ $ \varphi $ has an isolated real zero at the origin. Let $ s $ be a real number, with $ s\geq 1 $. Then we have the inclusion $ \imax_M^\infty \subseteq \varphi\mathcal{E}_n(M^s) $ if and only if $ s\geq \tau(\varphi) $.
\end{thm}

Together with Corollary \ref{ellip}, Theorem \ref{ellipt} yields a closedness result.

\begin{cor}\label{zis}
Assume that the sequence $ M $ is strongly regular and that the real-analytic germ $ \varphi $ has an isolated real zero at the origin. Then the ideal  $ \mathcal{I}=\varphi\mathcal{E}_n(M) $ is closed if and only if $ \tau(\varphi)=1 $.
\end{cor} 

\begin{exam}\label{exzerois}
Put $ n=2 $, $ (x_1,x_2)=(x,y) $. For a given integer $ k\geq 1 $, put $ \varphi(x,y)=y^2+x^{2k} $. Then $ \tau(\varphi)=k $, hence the ideal generated by $\varphi $ in a strongly regular local ring $ \mathcal{E}_2(M) $ is closed if and only if $ k=1 $. 
\end{exam}

\begin{exam}
The ``only if'' part of Theorem \ref{ellipt} can be illustrated as follows. As in Example \ref{exzerois}, put $ n=2 $, $ (x_1,x_2)=(x,y) $ and $ \varphi(x,y)=y^2+x^{2k} $. Consider a Gevrey sequence $ M_j=(j!)^\alpha $ with $ \alpha>0 $ and let $ h $ be the smooth function defined on $ \mathbb{R}^2 $ by $ h(x,y)=\exp(-\vert x\vert^{-1/\alpha}) $ for $ x\neq 0 $ and $ h(0,y)=0 $. The fact that $ h $ defines an element of $\imax^\infty_M $ is classical\footnote{In order to majorize the derivatives of $ h  $ at $ (x,y) $, with $ x\neq 0 $, one uses the Cauchy formula on a disk centered at $ x$ with radius $ \delta \vert x\vert $ for some sufficiently small $ \delta $ (depending on $\alpha $). This yields a bound of the form $ A^jj!\vert x\vert^{-j}\exp(-A/\vert x\vert^{-1/\alpha}) $ for the derivatives of order $ j $. It is then enough to write $ \vert x\vert^{-j}=B^j(j!)^\alpha\big((B\vert x\vert)^{-j/\alpha}/j!\big)^\alpha \leq B^j(j!)^\alpha\exp\big(\alpha(B\vert x\vert)^{-1/\alpha}\big) $ with $ B=(\alpha/A)^\alpha $.}. We claim, however, that $ h $ does not belong to $ \varphi\mathcal{E}_2(M^s) $ for $ s<k $. Indeed, consider the smooth function $ g=h/\varphi $. For $ 0\leq \vert y\vert <x^{2k} $, we have the expansion
$ g(x,y)=\sum_{j= 0}^{+\infty}(-1)^j\exp(-\vert x\vert^{-1/\alpha})x^{-2k(j+1)}y^{2j} $, hence
$ {\partial^{2j} g}/{\partial y^{2j}}(x,0)=(-1)^j(2j)!\exp(-\vert x\vert^{-1/\alpha})x^{-2k(j+1)} $. Using Stirling's formula, we derive
\begin{equation*}
\vert{\partial^{2j} g}/{\partial y^{2j}}(j^{-\alpha},0)\vert\geq C^{j+1}(2j)!^{1+k\alpha}=C^{j+1}(2j)!(M_{2j})^k
\end{equation*}
for some suitable constant $ C>0 $ and for any integer $ j\geq 1 $. The claim readily follows.
\end{exam}

Corollary \ref{zis} involves a rather restrictive assumption on the real zero set of $ \varphi $. If no such assumption is made, it is not clear whether the closedness of $ \varphi\mathcal{E}_n(M) $ can be analogously characterized by simple geometric data. However, the two-dimensional case is completely understood: we describe the solution hereafter.

Denote by $ \omega(v) $ the order of any non-zero element $ v $ of $ \mathcal{O}_1 $, that is, the smallest degree of the monomials in the power series expansion of $ v $. We put $ \omega(0)=\infty $. Now, put $ n=2 $, $ (x_1,x_2)=(x,y) $, and consider a non-zero element $ \varphi $ of $ \mathcal{O}_2 $. After a real linear change of coordinates, we can assume that the complex zero set of $ \varphi $ is not tangent to the $ y $-axis. Using Puiseux's theorem, we can then write 
\begin{equation*}
\varphi(x^m,y)=u(x^m,y)\prod_{j=1}^{k}\big(y-y_j(x)\big)^{n_j}, 
\end{equation*}
where $ m $ and the $ n_j $'s are integers with $ m\geq 1 $ and $ n_j\geq 1 $, the germ $ u $ is a unit in $ \mathcal{O}_2 $ and the $ y_j $'s are elements of $ \mathcal{O}_1 $, satisfying $ \omega(y_j)\geq m $ for $ 1\leq j\leq k $. With any root $ y_j $ we associate an integer $ d_j $ defined as follows: $ d_j=1 $ if the imaginary part $ \Im y_j $ is identically zero, $ d_j=\omega(\Im y_j) $ otherwise. We now put $ d^+(\varphi)=\max_{1\leq j\leq k} d_j/m $ and $ d(\varphi)=\max\big(d^+(\varphi),d^+(\check\varphi)\big) $, 
where $ \check\varphi(x,y)=\varphi(-x,y) $. 

\begin{exam}
For $ \varphi(x,y)=y^2-x^4 $, we have $ d(\varphi)=1 $, whereas for $ \varphi(x,y)=y^2+x^4 $, we have $ d(\varphi)=2 $. For $ \varphi(x,y)=y^2-x^3 $, we have $ d(\varphi)=3/2 $, as for $ \varphi(x,y)=y^2+x^3 $.
\end{exam}

Denote by $ \mathcal{I}' $ the set of those elements of $ \mathcal{E}_2(M) $ which are the germ of an element of $ \overline{\mathcal{I}(\nu)} $ for some $ \nu\geq 1 $. The following result can be found in \cite{VT2}.

\begin{thm}
Assume that $ n=2 $ and that the sequence $ M $ is strongly regular. Let $ s $ be a real number, with $ s\geq 1 $. Then we have $ \mathcal{I}'\subseteq \varphi\mathcal{E}_2(M^s) $ if and only if $ s\geq d(\varphi) $.
\end{thm}

The theorem has an obvious corollary.

\begin{cor}\label{n2}
Assume that $n=2 $ and that the sequence $ M $ is strongly regular. Then the ideal  $ \mathcal{I}=\varphi\mathcal{E}_2(M) $ is closed if and only if $ d(\varphi)=1 $.
\end{cor}
 
\begin{rem}
When $ \varphi $ has an isolated real zero at $ 0 $ in $ \mathbb{R}^2 $, the number $ d(\varphi) $ coincides with the \L ojasiewicz exponent $ \tau(\varphi) $. Thus, in this particular case, Corollary \ref{n2} agrees (as one could hope!) with Corollary \ref{zis}. But for general zero sets, we only have $ \tau(\varphi) \leq d(\varphi) $, and the inequality can be strict: see \cite{VT2}, Example 1.5 (iv). 
\end{rem}

\section{A glimpse on other results} 

After the strongly non-quasianalytic considerations of \S \ref{off}, let us come back to the quasianalytic case. Although we cannot give here a complete picture of the subject, it is nonetheless necessary to say a few words about an important recent addition to the theory, namely \emph{resolution of singularities}, for which we refer the reader to the articles of Bierstone and  Milman \cite{Bierstone-Milman} or Rolin, Speissegger and Wilkie \cite{RSW}. Once such a powerful tool is available, many applications arise.
For instance, although we do not know whether quasianalytic rings $ \mathcal{E}_n(M) $ are noetherian, they satisfy a weaker property of \emph{topological noetherianity}. This property involves $ M $-quasianalytic sets, whose definition is a straightforward extension of the analytic case\footnote{For any point $ a $ of $ \mathbb{R}^n $, put $ {}_a\mathcal{E}_n(M)=\{f(x-a) : f\in\mathcal{E}_n(M)\} $. A subset $ X $ of $ \mathbb{R}^n $ is called $ M$-quasianalytic if, for any $ a \in \mathbb{R}^n $, there is a finite family $ f_1,\ldots,f_p $ of elements of ${}_a\mathcal{E}_n(M) $ and an open neighborhood $ U $ of $ a $ such that $ X\cap U=\{x\in U:f_1(x)=\cdots=f_p(x)=0\} $.}. The result can be  stated as follows. 

\begin{thm}
Assume that the quasianalytic local ring $ \mathcal{E}_n(M) $ is stable under derivation. Then any decreasing sequence $ X_0\supseteq X_1\supseteq X_2\supseteq\cdots $ of $ M $-quasianalytic subsets of $ \mathbb{R}^n $ stabilizes in some neighborhood of a given compact set. 
\end{thm}

Resolution of singularities also provides a curve selecting lemma and, just as in the analytic case (see, for instance, \cite{Bochnak-Risler}), \L ojasiewicz inequalities can then be obtained as an application of curve selection.  

\begin{thm}
Assume that the quasianalytic local ring $ \mathcal{E}_n(M) $ is stable under derivation, and 
let $ f $ and $ g $ be two elements of $ \mathcal{E}_n(M) $ with $ g^{-1}(\{0\})\subseteq f^{-1}(\{0\}) $. Then there exist positive constants $ C>0 $ and $ \alpha \geq 1 $ such that $ \vert g(x)\vert\geq C\vert f(x)\vert^\alpha $ for any $ x $ in a neighborhood of $ 0 $.
\end{thm} 

Among other recent developments on quasianalyticity from the geometric viewpoint, we mention, finally, interesting connections with potential theory in several variables \cite{Coman1,Coman2,Pierzchala,Plesniak}, some of them being closely related to the aforementioned desingularization results.


\begin{thebibliography}{00}

\bibitem{Bang} \textsc{T. Bang}, \textit{Om quasi-analytiske Funktioner}, Thesis, University of Copenhagen, 1946.

\bibitem{Beaugendre} \textsc{P. Beaugendre}, \textit{Intersections de classes non-quasianalytiques}, Thesis, University of Paris XI-Orsay, 2002. Available online: 

\texttt{http://tel.ccsd.cnrs.fr/docs/00/04/48/43/PDF/tel-00001335.pdf}.

\bibitem{Bierstone-Milman} \textsc{E. Bierstone \& P.D. Milman}, \textit{Resolution of singularities in Denjoy-Carleman classes}, Selecta Math. \textbf{10} (2004), 1--28.

\bibitem{Bochnak-Risler} \textsc{J. Bochnak \& J.-J. Risler}, \textit{Sur les exposants de \L ojasiewicz}, Comment. Math. Helvetici \textbf{50} (1975), 493--507.

\bibitem{BBMT} \textsc{J. Bonet, R.W. Braun, R. Meise \& B.A. Taylor}, \textit{Whitney's extension theorem for nonquasianalytic classes of ultradifferentiable functions}, Studia Math. \textbf{99} (1991), 155--184.

\bibitem{Borel1} \textsc{\'E. Borel}, \textit{Sur la g\'en\'eralisation du prolongement analytique}, C. R. Acad. Sci. Paris \textbf{130} (1900), 1115--1118.

\bibitem{Borel2} \textsc{\'E. Borel}, \textit{Sur les s\'eries de polyn\^omes et de fractions rationnelles}, Acta Math. \textbf{24} (1901), 309--387.

\bibitem{Bourbaki} \textsc{N. Bourbaki}, \textit{Espaces Vectoriels Topologiques, Chapitres 1 \`a 5}, Masson, 1981.

\bibitem{Carleman1} \textsc{T. Carleman}, \textit{Sur les fonctions quasi-analytiques},  Comptes Rendus du 5\`eme congr\`es des math\'ema\-ticiens scandinaves \`a Helsingfors (1922), 181--196. 

\bibitem{Carleman11} \textsc{T. Carleman}, \textit{Sur le calcul effectif d'une fonction quasi-analytique dont on donne les d\'eriv\'ees en un point}, C. R. Acad. Sci. Paris \textbf{176} (1923), 64--65.

\bibitem{Carleman12} \textsc{T. Carleman}, \textit{Sur les fonctions ind\'efiniment d\'erivables}, C. R. Acad. Sci. Paris \textbf{177} (1923), 422--424.


\bibitem{Carleman2} \textsc{T. Carleman}, \textit{Les fonctions quasi-analytiques}, Gauthiers Villars, Paris, 1926.

\bibitem{Cartan1} \textsc{H. Cartan}, \textit{Sur les classes de fonctions d\'efinies par des in\'egalit\'es portant sur leurs d\'eriv\'ees successives}, Actual. Sci. Ind. \textbf{867}, Hermann, Paris, 1940.

\bibitem{Cartan2} \textsc{H. Cartan, S. Mandelbrojt}, \textit{Solution du probl\`eme d'\'equivalence des classes de fonctions ind\'efiniment d\'erivables}, 
Acta Math. \textbf{72} (1940), 31--49.

\bibitem{Childress} \textsc{C.L. Childress}, \textit{Weierstrass division in quasianalytic local rings}, Canad. J. Math. \textbf{28} (1976), 938--953.

\bibitem{ChtCht0} \textsc{J. Chaumat, A.-M. Chollet}, \textit{Sur le th\'eor\`eme de division de Weierstrass}, Studia Math. \textbf{116} (1995), 59--84.

\bibitem{ChtCht1} \textsc{J. Chaumat, A.-M. Chollet}, \textit{Caract\'erisation des anneaux noetheriens de s\'eries formelles \`a croissance contr\^ol\'ee. Application \`a la synth\`ese spectrale}, Publ. Mat. \textbf{41} (1997), 545--561.

\bibitem{ChtCht} \textsc{J. Chaumat, A.-M. Chollet}, \textit{Division par un polyn\^ome hyperbolique}, Canad. J. Math. \textbf{56} (2004), 1121--1144.

\bibitem{Coman1} \textsc{D. Coman, N. Levenberg, E. Poletsky}, \textit{Quasianalyticity and pluripolarity}, J. Amer. Math. Soc. \textbf{18} (2005), 239--252.

\bibitem{Coman2} \textsc{D. Coman, N. Levenberg, E. Poletsky}, \textit{Smooth submanifolds intersecting any analytic curve in a discrete set}, Math. Ann. \textbf{332} (2005), 55--65.

\bibitem{Denjoy} \textsc{A. Denjoy}, \textit{Sur les fonctions quasi-analytiques de variable r\'eelle}, C. R. Acad. Sci. Paris \textbf{123} (1921), 1320--1322.

\bibitem{Grauert} \textsc{H. Grauert, R. Remmert}, \textit{Analytische Stellenalgebren}, Springer-Verlag, Berlin, 1971.

\bibitem{Gevrey} \textsc{M. Gevrey}, \textit{Sur la nature analytique des solutions des \'equations aux d\'eriv\'ees partielles}, Ann. Sci. \'Ec. Norm. Sup. \textbf{35} (1918), 129--190.

\bibitem{Hadamard} \textsc{J. Hadamard}, \textit{Sur la g\'en\'eralisation de la notion de fonction analytique}, Bull. Soc. Math. France \textbf{40} (1912), suppl\'ement sp\'ecial~: vie de la soci\'et\'e, s\'eance du 28 f\'evrier 1912, 28-29.

\bibitem{Hormander} \textsc{L. H\"ormander}, \textit{The Analysis of Linear Partial Differential Operators I}, Springer Verlag, Berlin, 1990.

\bibitem{Katznelson} \textsc{Y. Katznelson}, \textit{An introduction to harmonic analysis}, Dover, New York, 1976.

\bibitem{Komatsu} \textsc{H. Komatsu}, \textit{Ultradistributions, I. Structure theorems and a characterization}, J. Fac. Sci. Univ. Tokyo, Sect. IA \textbf{20} (1973), 25--105.

\bibitem{Komatsu2} \textsc{H. Komatsu}, \textit{The implicit function theorem for ultradifferentiable mappings}, Proc. Japan Acad. ser. A \textbf{55} (1979), 69--72.

\bibitem{Koosis} \textsc{P. Koosis}, \textit{The logarithmic integral, I}, Cambridge University Press, Cambridge, 1998.

\bibitem{Malgrange} \textsc{B. Malgrange}, \textit{Ideals of Differentiable Functions}, Tata Institute of Fundamental Research, Bombay; Oxford University Press, London, 1967.

\bibitem{Malgrange2} \textsc{B. Malgrange}, \textit{Id\'eaux de fonctions diff\'erentiables et division des distributions}, in: \textit{Dans le sillage de Laurent Schwartz}, Actes des journ\'ees X-UPS 2003, \'Ecole Polytechnique, Palaiseau, 2003. Available online at:

\texttt{http://www.math.polytechnique.fr/xups/vol03.html}

\bibitem{Mandelbrojt} \textsc{S. Mandelbrojt}, \textit{Sur les fonctions ind\'efiniment d\'erivables}, Acta Math. \textbf{72} (1940), 15--29.

\bibitem{Meise-Vogt} \textsc{R. Meise, D. Vogt}, \textit{Introduction to Functional Analysis}, The Clarendon Press, Oxford University Press, New York, 1997.

\bibitem{Mouze} \textsc{A. Mouze}, \textit{Anneaux de s\'eries formelles \`a croissance contr\^ol\'ee}, Thesis, University of Lille, 2000. Available online at:

\texttt{http://tel.ccsd.cnrs.fr/docs/00/08/03/23/PDF/TheseAM.pdf}

\bibitem{Petzsche} \textsc{H.-J. Petzsche}, \textit{On E. Borel's theorem}, Math. Ann. \textbf{282} (1988), 299--313.

\bibitem{Pierzchala} \textsc{R. Pierzcha\l a}, \textit{UPC condition in polynomially bounded o-minimal structures}, J. Approx. Theory. \textbf{132} (2005), 25-33.

\bibitem{Plesniak} \textsc{W. Ple\'sniak}, \textit{Pluriregularity in polynomially bounded o-minimal structures}, Univ. Iagel. Acta Math. \textbf{41} (2003), 205--214.

\bibitem{Remmert} \textsc{R. Remmert}, \textit{Classical topics in complex function theory}, Springer Verlag, Berlin, 1998.

\bibitem{RSW} \textsc{J.-P. Rolin, P. Speissegger, A. Wilkie}, \textit{Quasianalytic Denjoy-Carleman classes and o-minimality}, J. Amer. Math. Soc. \textbf{16} (2003), 751--777.

\bibitem{Roumieu} \textsc{C. Roumieu}, \textit{Ultradistributions d\'efinies sur $ \mathbb{R}^n $ et sur certaines classes de vari\'et\'es dif\-f\'e\-ren\-tia\-bles}, J. Analyse Math. \textbf{10} (1962-1963), 153--192.

\bibitem{Rudin} \textsc{W. Rudin}, \textit{Real and Complex Analysis}, MacGraw-Hill, New-York, 1987.

\bibitem{Schmets} \textsc{J. Schmets, M. Valdivia}, \textit{On certain extension theorems in the mixed Borel setting}, J. Math. Anal. Appl. \textbf{297} (2004), 384--403.

\bibitem{VT1} \textsc{V. Thilliez}, \textit{On closed ideals in smooth classes}, Math. Nachr. \textbf{227} (2001), 143--157.

\bibitem{VT2} \textsc{V. Thilliez}, \textit{A sharp division estimate for ultradifferentiable germs}, Pacific J. Math. \textbf{205} (2002), 237--256.

\bibitem{VT3} \textsc{V. Thilliez}, \textit{Bounds for quotients in rings of formal power series with growth constraints}, Studia Math. \textbf{151} (2002), 49--65.

\bibitem{Tougeron} \textsc{J.-C. Tougeron}, \textit{Id\'eaux de fonctions diff\'erentiables}, Springer Verlag, Berlin, 1972.

\end{thebibliography}
\end{document}